# Truncated transparent boundary conditions

Ivan Sofronov[1,2]

[1]*Schlumberger, Moscow, Pudovkina 13,*

[2]*MIPT, Moscow region, Dolgoprudny, Institutskii per. 1*

isofronov@slb.com

**Abstract**

We derive equations for evaluating differential operators in transparent boundary conditions (TBCs) for a certain class of hyperbolic systems of second-order equations. This local part of TBCs can be used as approximate nonreflecting boundary conditions. We give examples of obtaining such truncated TBCs for 3D elasticity and Biot poroelasticity.

*Keywords:* Transparent boundary conditions, local boundary conditions, hyperbolic systems, cylindrical anisotropy, orthotropic elasticity, Biot poroelasticity,

# 1   Introduction

Transparent boundary conditions (TBCs) for the wave equation were proposed in [10], [5] and were further developed, for example, in [12], [6], [1], [17]. In particular, on a circle ($d=2$) or a sphere ($d=3$) they can be written as

$$\frac{\partial u}{\partial t} + \frac{\partial u}{\partial r} + \frac{d-1}{2}\frac{c}{r}u - Q^{-1}\{B_k *\}Qu = 0 \qquad (1)$$

where $Q$ and $Q^{-1}$ are the Fourier and inverse Fourier transform operators over the trigonometric ($d=2$) or spherical ($d=3$) basis, and $\{B_k *\}$ denotes the time convolution operator of the $k$-th Fourier coefficient with kernels $B_k(t)$ derived analytically [12]. In order to provide efficient calculations with reasonable computational recourses the kernels $B_k(t)$ are approximated by sums of exponentials,

$$B_k(t) \approx \sum_{l=1}^{L_k} a_{l,k} e^{b_{l,k}t}, \quad \operatorname{Re} b_{l,k} \leq 0 .$$

The explicit analytical expressions for other examples of TBCs operators can be obtained only for simple boundaries (plane, sphere, etc.) and for definite narrow classes of equations. Extension to the equations with variable coefficients is possible by using so-called *quasi-analytic* TBCs proposed in [13], [15] with numerically found operator matrices.

The pseudo-differential TBCs operators consist of two parts: a local (differential) part and a nonlocal part in space and time. For instance, in l.h.s. (1) the last term is computed by the non-local in both space and time operator, while the three first terms define action of a local operator. In [14] there were announced equations for evaluating this local part of TBC operator for general second-order hyperbolic systems. This is a natural idea to use it as an approximate operator of local non-reflecting boundary conditions while solving external initial boundary-valued problems. We call such condition *truncated* TBC (TTBC). Note that, e.g., for the wave equation TTBCs coincide with the well-known characteristic boundary conditions. Practical application of TTBCs derived in [14] for 3D elasticity equations is studied in [16].

In this paper, Section 1, we present the derivation of the equations announced in [14] and consider the way of their solving, i.e. to evaluate operators of the TBCs' differential part.

In Section 2, we consider another two examples of generating TTBCs: for 3D orthotropic elasticity equations in cylindrical coordinates, and for 3D Biot poroelasticity equations in Cartesian coordinates.

## 2   Derivation of equations for operators of TTBCs

### *2.1   Problem formulation*

Consider a domain $\Omega \subset R^3$ with a sufficiently smooth open boundary $\Gamma$ at the point $O$. Denote by $n$ the outward normal and by $\tau = (\tau_1, \tau_2)$ the vector of orthogonal coordinates in a tangent plane, see Fig.1.

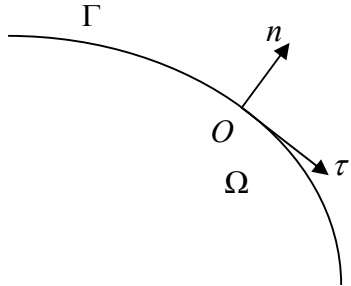

Fig. 1. External boundary $\Gamma$, local coordinates at the point $O$

Let a vector function $u = (u^1, u^2, ..., u^N)^\mathrm{T}$ satisfy a second-order hyperbolic system of equations [4] written in local coordinates $(n, \tau)$ at the point $O$ in the form

$$-\frac{\partial^2 u}{\partial t^2} + A\frac{\partial^2 u}{\partial n^2} + (B\nabla_\tau)\frac{\partial u}{\partial n} + C\frac{\partial u}{\partial n} + a(\nabla_\tau, \nabla_\tau)u + (c\nabla_\tau)u + du = 0 . \quad (2)$$

Here, we use the following notation:

$$(B\nabla_\tau)u \equiv \sum_{i=1}^{2} B_i \frac{\partial u}{\partial \tau_i}, \quad a(\nabla_\tau, \nabla_\tau)u \equiv \sum_{i=1}^{2}\sum_{j=1}^{2} a_{ij}\frac{\partial^2 u}{\partial \tau_i \partial \tau_j}.$$

The $N \times N$ matrix coefficients $A(\tau)$, $B_i(\tau)$ and $C(n,\tau)$, $a_{ij}(n,\tau)$, $c_i(n,\tau)$, and $d(n,\tau)$ may depend on $\tau$ and on $n$, $\tau$, respectively.

We suppose that the coefficients and solution of (2) are sufficiently smooth in the neighborhood of the point $O$.

Our goal is to obtain a relationship between the solution $u$ and its first derivatives, which is approximately valid at the point $O$. This relationship will have the form $P_1 \frac{\partial u}{\partial t} - \frac{\partial u}{\partial n} + P_0 u = 0$ with some operators $P_1$, $P_0$ that are obtained from the coefficients of (2). Initial condition of rest at infinity will be taken into account in order to discriminate the wave solutions $u$ leaving the domain $\Omega$.

## 2.2 Derivation of equations for $P_1$, $P_0$

Applying the Laplace transform with respect to time to (2) we have

$$-s^2\hat{u} + A\frac{\partial^2 \hat{u}}{\partial n^2} + (B\nabla_\tau)\frac{\partial \hat{u}}{\partial n} + C\frac{\partial \hat{u}}{\partial n} + a(\nabla_\tau, \nabla_\tau)\hat{u} + (c\nabla_\tau)\hat{u} + d\hat{u} = 0 .$$

After introducing the vanishing parameter $\varepsilon = 1/s \to +0$ and local coordinate transformation $\eta = n/\varepsilon$ we obtain for the function $\tilde{u}(\eta,\tau) \equiv \hat{u}(\varepsilon\eta,\tau)$ the equation

$$A\frac{\partial^2 \tilde{u}}{\partial \eta^2} + \varepsilon(B\nabla_\tau)\frac{\partial \tilde{u}}{\partial \eta} + \varepsilon\tilde{C}\frac{\partial \tilde{u}}{\partial \eta} + \varepsilon^2 \tilde{a}(\nabla_\tau, \nabla_\tau)\tilde{u} + \varepsilon^2(\tilde{c}\nabla_\tau)\tilde{u} + \varepsilon^2\tilde{d}\tilde{u} - \tilde{u} = 0 \quad (3)$$

(the tilde in the coefficients denotes that they are functions of $\eta$).

Since $\frac{\partial}{\partial n}C(n,\tau)$ is bounded in physical coordinates $(n,\tau)$, we evidently have as $\varepsilon \to 0$:

$$\frac{\partial}{\partial \eta}\tilde{C}(\eta,\tau)|_{\eta=0} = \varepsilon \frac{\partial}{\partial n}C(n,\tau)|_{n=0} = O(\varepsilon). \tag{4}$$

Denote by $P_1 + \varepsilon P_0$ the approximate Poincare-Steklov operator for (3) written as a linear function of $\varepsilon$, i.e., consider the following relationship between $\tilde{u}$ and its normal derivative:

$$\frac{\partial}{\partial \eta}\tilde{u}(\eta,\tau,\varepsilon)|_{\eta=0} = (P_1 + \varepsilon P_0)\tilde{u}(\eta,\tau,\varepsilon)|_{\eta=0} + o(\varepsilon). \tag{5}$$

To obtain $P_1$ and $P_0$ we introduce the linear expansion of $\tilde{u}$ with respect to $\varepsilon$:

$$\tilde{u}(\eta,\tau,\varepsilon) = \tilde{v}(\eta,\tau) + \varepsilon \tilde{w}(\eta,\tau) + o(\varepsilon)$$

and substitute it for (3). Combining terms of zero and first degrees of $\varepsilon$ we obtain:

$$A\frac{\partial^2 \tilde{v}}{\partial \eta^2} - \tilde{v} + \varepsilon A \frac{\partial^2 \tilde{w}}{\partial \eta^2} - \varepsilon \tilde{w} + \varepsilon (B\nabla_\tau)\frac{\partial \tilde{v}}{\partial \eta} + \varepsilon C_0 \frac{\partial \tilde{v}}{\partial \eta} + o(\varepsilon) = 0. \tag{6}$$

Here we use $C_0(\tau) \equiv \tilde{C}(0,\tau)$ instead of $C(n,\tau)$ according to (4).

Taking the limit as $\varepsilon \to 0$ in (6), we come to the following necessary conditions:

$$A\frac{\partial^2 \tilde{v}}{\partial \eta^2} - \tilde{v} = 0 \tag{7}$$

and

$$A\frac{\partial^2 \tilde{w}}{\partial \eta^2} - \tilde{w} + \left((B\nabla_\tau) + C_0\right)\frac{\partial \tilde{v}}{\partial \eta} = 0. \tag{8}$$

The waves leaving $\Omega$ correspond to vanishing branch of solutions to (7) as $\eta \to +\infty$. These solutions have the form

$$\tilde{v}(\eta,\tau) = \exp\left(-A^{-1/2}(\tau)\eta\right)\tilde{v}_0(\tau)$$

where $A^{-1/2}(\tau)$ is a positive square root of $A^{-1}(\tau)$, and $\tilde{v}_0(\tau)$ is a function.

Thus we obtain the following Poincare-Steklov map for $\tilde{v}$:

$$\frac{\partial \tilde{v}(\eta,\tau)}{\partial \eta}\bigg|_{\eta=0} = -A^{-1/2}(\tau)\tilde{v}(\eta,\tau)\big|_{\eta=0} . \qquad (9)$$

Consider now (8) as equation with respect to $\tilde{w}$ supposing that $\dfrac{\partial \tilde{v}}{\partial \eta}$ is known. It is the non-uniform (7) with a specific r.h.s. Denoting for brevity

$$E(\eta,\tau) \equiv \exp\left(-A^{-1/2}(\tau)\eta\right)$$

we transform the r.h.s. as follows:

$$\begin{aligned}
\left((B\nabla_\tau)+C_0\right)\frac{\partial \tilde{v}}{\partial \eta} &= \left((B\nabla_\tau)+C_0\right)\frac{\partial}{\partial \eta} E(\eta,\tau)\tilde{v}_0(\tau) \\
&= -(B\nabla_\tau)A^{-1/2}(\tau)E(\eta,\tau)\tilde{v}_0(\tau) - C_0 A^{-1/2}(\tau)E(\eta,\tau)\tilde{v}_0(\tau) \\
&= -(B\nabla_\tau A^{-1/2})E\tilde{v}_0 - (BA^{-1/2}\nabla_\tau E)\tilde{v}_0 - (BA^{-1/2}E\nabla_\tau)\tilde{v}_0 - C_0 A^{-1/2}E\tilde{v}_0 \\
&= -\left((B\nabla_\tau A^{-1/2})E + (BA^{-1/2}E\nabla_\tau) + C_0 A^{-1/2}E\right)\tilde{v}_0 + O(\eta)
\end{aligned}$$

Evidently $E(0,\tau) = I$; therefore, (8) at $\eta = 0$ reads:

$$A\frac{\partial^2 \tilde{w}}{\partial \eta^2} - \tilde{w} - D\tilde{v}_0 = 0 \qquad (10)$$

where

$$D \equiv (B\nabla_\tau A^{-1/2}) + (BA^{-1/2}\nabla_\tau) + C_0 A^{-1/2} .$$

It follows from (7) and (10) that

$$\begin{aligned}
\frac{\partial^2}{\partial \eta^2}\tilde{u}(\eta,\tau,\varepsilon)\big|_{\eta=0} &= A^{-1}\tilde{v} + \varepsilon A^{-1}\left(\tilde{w} + D\tilde{v}_0\right) \\
&= A^{-1}\tilde{u}_0 + \varepsilon A^{-1}D\tilde{v}_0 \qquad (11)\\
&= A^{-1}\tilde{u}_0 + \varepsilon A^{-1}D\tilde{u}_0 + o(\varepsilon)
\end{aligned}$$

where $\tilde{u}_0 \equiv \tilde{u}(\eta, \tau, \varepsilon)|_{\eta=0}$.

Now notice that (6) is also valid for the normal derivative $\dfrac{\partial}{\partial \eta}\tilde{u}(\eta, \tau, \varepsilon)$ since the coefficients of (6) don't depend on $\eta$. Therefore, it follows from (5) that

$$\begin{aligned}\dfrac{\partial^2}{\partial \eta^2}\tilde{u}(\eta, \tau, \varepsilon)|_{\eta=0} &= (P_1 + \varepsilon P_0)\dfrac{\partial}{\partial \eta}\tilde{u}(\eta, \tau, \varepsilon)|_{\eta=0} + o(\varepsilon) \\ &= (P_1 + \varepsilon P_0)^2 \tilde{u}(\eta, \tau, \varepsilon)|_{\eta=0} + o(\varepsilon) \\ &= (P_1 P_1 + \varepsilon P_1 P_0 + \varepsilon P_0 P_1)\tilde{u}(\eta, \tau, \varepsilon)|_{\eta=0} + o(\varepsilon)\end{aligned} \qquad (12)$$

Using (11) and (12), we obtain

$$(P_1 P_1 + \varepsilon P_1 P_0 + \varepsilon P_0 P_1)\tilde{u}_0 = A^{-1}\tilde{u}_0 + \varepsilon A^{-1}D\tilde{u}_0 + o(\varepsilon) .$$

Equating terms of zero and first degrees of $\varepsilon$, we derive the simple equation for finding $P_1$:

$$P_1 P_1 = A^{-1} ,$$

and the Sylvester equation for finding $P_0$:

$$P_1 P_0 + P_0 P_1 = A^{-1}D .$$

Hence

$$P_1 = -A^{-1/2} , \qquad (13)$$

the negative sign follows from (9), and

$$A^{-1/2}P_0 + P_0 A^{-1/2} = -A^{-1}\left((B\nabla_\tau A^{-1/2}) + (BA^{-1/2}\nabla_\tau) + C_0 A^{-1/2}\right) . \qquad (14)$$

Coming back to the local coordinates $(n, \tau)$ in (5) we get as $s \to \infty$:

$$\dfrac{\partial \hat{u}}{\partial n}\bigg|_{n=0} = (sP_1 + P_0)\hat{u}|_{n=0} + o(s^0) . \qquad (15)$$

(we use here notation $1 = s^0$ to keep the variable in the asymptotically negligible expansion term).

Neglecting the vanishing term we formulate the relationship on the boundary:

$$\frac{\partial \hat{u}}{\partial n} = sP_1\hat{u} + P_0\hat{u} \ . \qquad (16)$$

After the inverse Laplace transform of (16) we obtain the following approximate relationship of solution to (2) and its derivatives at the point $O$ :

$$P_1 \frac{\partial u}{\partial t} - \frac{\partial u}{\partial n} + P_0 u = 0 \qquad (17)$$

with operators $P_1$ , $P_0$ , defined by (13) and the Sylvester equation (14).

Existence of operators $P_1$ and $P_0$ in (17) follows from the hyperbolicity of (2). Indeed, introducing the characteristic determinant

$$Q(\lambda, \xi_n, \xi_1, \xi_2) = \det\left(-\lambda^2 I + A\xi_n^2 + \sum_{i=1}^{2} B_i \xi_i \xi_n + \sum_{i=1}^{2} \sum_{j=1}^{2} a_{ij} \xi_i \xi_j \right),$$

see [4], [3], where $I$ is the unity $N \times N$ matrix, we have the hyperbolicity condition that the equation $Q(\lambda, \xi_n, \xi_1, \xi_2) = 0$ should have $2N$ different real roots $\lambda$ at the point $O$ for arbitrary real $\xi_n, \xi_1, \xi_2$ , $|\xi| = 1$ . Suppose $\xi_n = 1$ , $\xi_1 = \xi_2 = 0$ . Then existence of $2N$ different real roots $\lambda$ in the equation

$$\det\left(A - \lambda^2 I\right) = 0$$

means that all eigenvalues of matrix $A$ are positive [9]. I.e., it has the reciprocal squared root $A^{-1/2}$ . From positiveness of eigenvalues of matrix $A^{-1/2}$ it also follows the resolution of the Sylvester equation (14) with respect to $P_0$ [8].

## 2.3 Truncated TBCs

Evidently, a further more accurate treatment of the last term in (15) will generate non-local in time operators after the inverse Laplace transform. Thus the local in time part of the transparent (i.e., exact) boundary conditions is completely described in (17). Also note that all known operators of TBCs for the second order hyperbolic equations [10], [11], [12], [5], [6], [7] are the sum of the operator in (17) and a nonlocal operator (cf. (1)). The quasi-analytic TBCs, see details and notation in [13], [15], have a more complex but similar form to analytic TBC,

$$\mathbf{Q}^{-1}\mathbf{P_1Q}\frac{\partial f}{\partial t} - \frac{\partial f}{\partial r} + \mathbf{Q}^{-1}\mathbf{P_0Q}f + \mathbf{Q}^{-1}\left\{\tilde{\mathbf{K}}(\mathbf{t})*\right\}\mathbf{Q}f = 0 \ .$$

Equation (17) shows that one can use in this condition operators $P_1$ and $P_0$ instead of $\mathbf{Q}^{-1}\mathbf{P_1Q}$ and $\mathbf{Q}^{-1}\mathbf{P_0Q}$, respectively.

Therefore we suggest the equation (17) as an approximate local non-reflecting boundary condition for system (2), and call it Truncated TBCs (TTBC).

The analysis made in this Section proves the following theorem

**Theorem.** *The differential part of the transparent boundary conditions for system (2) has the form*

$$P_1\frac{\partial u}{\partial t} - \frac{\partial u}{\partial n} + P_0 u$$

*where*

$$P_1 = -A^{-1/2}$$

*and $P_0$ is determined by the Sylvester equation*

$$A^{-1/2}P_0 + P_0 A^{-1/2} = -A^{-1}\left((B\nabla_\tau A^{-1/2}) + (BA^{-1/2}\nabla_\tau) + C_0 A^{-1/2}\right) \ .$$

## 2.4 Evaluation of $P_1$, $P_0$

Let us describe details of computations of matrices of TTBC operator.

We are looking for operators $P_1$ and $P_0$ at the point $O$. The considered coefficients $A$, $B$, and $C$ in (2) can depend on $\tau$. Therefore, let us introduce the matrices with numerical entries:

$$A_O = A|_O, \; B_O = B|_O, \; C_O = C_0|_O \equiv C|_O \;.$$

Then (13) is replaced by

$$P_1 = -A_O^{-1/2} \;. \tag{18}$$

Thus the matrix $P_1$ is found by usual algorithms for positively defined matrices.

A more complex situation can arise for $P_0$ because it contains differential operators. Consider how to solve (14) for sufficiently smooth non-uniform coefficients.

We represent $P_0$ by sum of algebraic and differential terms:

$$P_0 = p + (q\nabla_\tau) \tag{19}$$

where matrices $p$ and $q = (q_1, q_2)$ will be separately calculated. Introduce for brevity the following matrices:

$$L = A_O^{-1/2}, \; M = -A_O^{-1} C_O A_O^{-1/2} - A_O^{-1}(B_O \nabla_\tau A_O^{-1/2}), \; K = -A_O^{-1}(B_O A_O^{-1/2}) = -A_O^{-1} B_O A_O^{-1/2},$$

i.e., $(K\nabla_\tau) = -A_O^{-1}(B_O A_O^{-1/2} \nabla_\tau)$. Then (14) takes the form

$$LP_0 + (P_0 A^{-1/2})\big|_O = M + (K\nabla_\tau) \;.$$

Using (19) we obtain

$$Lp + pL + L(q\nabla_\tau) + (q\nabla_\tau A^{-1/2})\big|_O = M + (K\nabla_\tau) \;.$$

After collecting terms with and without the gradient operator $\nabla_\tau$ we come to the equations:

$$L(q\nabla_\tau) + (qL\nabla_\tau) = (K\nabla_\tau)$$
$$Lp + pL = M - (q\nabla_\tau A^{-1/2})\big|_O$$

The first one gives the following algebraic equation for $q$:

$$Lq + qL = K,$$

and the second one – for $p$ (with already known $q$).

Coming back to the original notation we obtain two matrix equations for sequential resolving with respect to $p$ and $q$:

$$A_O^{-1/2} q + q A_O^{-1/2} = -A_O^{-1} B_O A_O^{-1/2} \tag{20}$$

$$A_O^{-1/2} p + p A_O^{-1/2} = -A_O^{-1} C_O A_O^{-1/2} - A_O^{-1}(B_O \nabla_\tau A^{-1/2})\big|_O - (q\nabla_\tau A^{-1/2})\big|_O \tag{21}$$

Formulas (18), (19) with the Silvester equations (20), (21) describe the whole process of generating our TTBC operator for (2).

Note that it is often convenient to use TTBC in the form resolved with respect to the time derivative:

$$\frac{\partial u}{\partial t} - P_1^{-1} \frac{\partial u}{\partial n} + P_1^{-1} P_0 u = 0$$

or, denoting $p_1 = P_1^{-1}$, $p_0 = P_1^{-1} P_0$,

$$\frac{\partial u}{\partial t} - p_1 \frac{\partial u}{\partial n} + p_0 u = 0. \tag{22}$$

Then we straightforwardly have from (13), (14) the following equations for direct computation of $p_1$, $p_0$:

$$p_1 = -A^{1/2},$$

$$A^{-1/2} p_0 + p_0 A^{-1/2} = A^{-1/2}\left((B\nabla_\tau A^{-1/2}) + (BA^{-1/2}\nabla_\tau) + C_0 A^{-1/2}\right).$$

# 3  Examples of generating TTBCs

In addition to TTBCs derived in [14] for 2D VTI cylindrical anisotropy and 3D Navier wave equation we consider another two examples for elasticity and poroelasticity.

## 3.1  Three-dimensional orthotropic elasticity in cylindrical coordinates

Here we derive TTBCs on the side cylinder boundary $r = const$ for constant elastic coefficients (formulas at the boundaries $z = const$ can be similarly derived).

Let $U = (u, v, w)^T$ be the displacements vector in the cylindrical system of coordinates $(r, \theta, z)$. The uniform elasticity equations have the following form:

$$\begin{cases} \rho \dfrac{\partial^2 u}{\partial t^2} = \dfrac{\partial \sigma_r}{\partial r} + \dfrac{1}{r}\dfrac{\partial \tau_{r\theta}}{\partial \theta} + \dfrac{\partial \tau_{rz}}{\partial z} + \dfrac{\sigma_r - \sigma_\theta}{r} \\ \rho \dfrac{\partial^2 v}{\partial t^2} = \dfrac{\partial \tau_{r\theta}}{\partial r} + \dfrac{1}{r}\dfrac{\partial \sigma_\theta}{\partial \theta} + \dfrac{\partial \tau_{\theta z}}{\partial z} + \dfrac{2\tau_{r\theta}}{r} \\ \rho \dfrac{\partial^2 w}{\partial t^2} = \dfrac{\partial \tau_{rz}}{\partial r} + \dfrac{1}{r}\dfrac{\partial \tau_{\theta z}}{\partial \theta} + \dfrac{\partial \sigma_z}{\partial z} + \dfrac{\tau_{rz}}{r} \end{cases}$$

where, according to the Hooke's law for cylindrically orthotropic materials,

$$\begin{aligned} \sigma_r &= A_{11}\varepsilon_r + A_{12}\varepsilon_\theta + A_{13}\varepsilon_z \\ \sigma_\theta &= A_{12}\varepsilon_r + A_{22}\varepsilon_\theta + A_{23}\varepsilon_z \\ \sigma_z &= A_{13}\varepsilon_r + A_{23}\varepsilon_\theta + A_{33}\varepsilon_z \\ \tau_{\theta z} &= A_{44}\gamma_{\theta z}, \; \tau_{rz} = A_{55}\gamma_{rz}, \; \tau_{r\theta} = A_{66}\gamma_{r\theta} \end{aligned},$$

$\rho$ is the density, and $A_{11} \div A_{66}$ are the stiffness coefficients; the strains are defined as usual

$$\begin{aligned} \varepsilon_r &= \dfrac{\partial u}{\partial r}, \; \varepsilon_\theta = \dfrac{1}{r}\dfrac{\partial v}{\partial \theta} + \dfrac{u}{r}, \; \varepsilon_z = \dfrac{\partial w}{\partial z} \\ \gamma_{r\theta} &= \dfrac{\partial v}{\partial r} - \dfrac{v}{r} + \dfrac{1}{r}\dfrac{\partial u}{\partial \theta}, \; \gamma_{rz} = \dfrac{\partial w}{\partial r} + \dfrac{\partial u}{\partial z}, \; \gamma_{\theta z} = \dfrac{\partial v}{\partial z} + \dfrac{1}{r}\dfrac{\partial w}{\partial \theta} \end{aligned}.$$

Following notation used in (2), we collect from above equations the terms containing derivatives with respect to $r$ in the second order system:

$$-\frac{\partial^2 U}{\partial t^2} + A\frac{\partial^2 U}{\partial r^2} + B_\theta \frac{\partial^2 U}{\partial \theta \partial r} + B_z \frac{\partial^2 U}{\partial z \partial r} + C\frac{\partial U}{\partial r} + \ldots = 0, \qquad (23)$$

i.e., obtain the equations

$$\begin{cases} -\rho\dfrac{\partial^2 u}{\partial t^2} + A_{11}\dfrac{\partial^2 u}{\partial r^2} + (A_{12} + A_{66})\dfrac{1}{r}\dfrac{\partial^2 v}{\partial r \partial \theta} + (A_{13} + A_{55})\dfrac{\partial^2 w}{\partial r \partial z} + A_{11}\dfrac{1}{r}\dfrac{\partial u}{\partial r} + \ldots = 0 \\ -\rho\dfrac{\partial^2 v}{\partial t^2} + A_{66}\dfrac{\partial^2 v}{\partial r^2} + A_{66}\dfrac{1}{r}\dfrac{\partial v}{\partial r} + (A_{12} + A_{66})\dfrac{1}{r}\dfrac{\partial^2 u}{\partial r \partial \theta} + \ldots = 0 \\ -\rho\dfrac{\partial^2 w}{\partial t^2} + A_{55}\dfrac{\partial^2 w}{\partial r^2} + A_{55}\dfrac{1}{r}\dfrac{\partial w}{\partial r} + (A_{13} + A_{55})\dfrac{\partial^2 u}{\partial r \partial z} + \ldots = 0 \end{cases}.$$

Thus the $3\times 3$ matrices in (23) are as follows

$$A = \frac{1}{\rho}\begin{pmatrix} A_{11} & 0 & 0 \\ 0 & A_{66} & 0 \\ 0 & 0 & A_{55} \end{pmatrix}, \quad C = \frac{1}{\rho r}\begin{pmatrix} A_{11} & 0 & 0 \\ 0 & A_{66} & 0 \\ 0 & 0 & A_{55} \end{pmatrix},$$

$$B_\theta = \frac{1}{\rho r}\begin{pmatrix} 0 & A_{12} + A_{66} & 0 \\ A_{12} + A_{66} & 0 & 0 \\ 0 & 0 & 0 \end{pmatrix}, \quad B_z = \frac{1}{\rho}\begin{pmatrix} 0 & 0 & A_{13} + A_{55} \\ 0 & 0 & 0 \\ A_{13} + A_{55} & 0 & 0 \end{pmatrix}.$$

Using (13), (14) we calculate the matrix

$$P_1 = -A^{-1/2} = -\sqrt{\rho}\begin{pmatrix} A_{11}^{-1/2} & 0 & 0 \\ 0 & A_{66}^{-1/2} & 0 \\ 0 & 0 & A_{55}^{-1/2} \end{pmatrix},$$

and obtain the equation for $P_0$ :

$$A^{-1/2}P_0 + P_0 A^{-1/2} = -\frac{1}{r}\rho^{-1}A^{-1}\begin{pmatrix} 0 & A_{12}+A_{66} & 0 \\ A_{12}+A_{66} & 0 & 0 \\ 0 & 0 & 0 \end{pmatrix} A^{-1/2}\frac{\partial}{\partial \theta}$$

$$-\rho^{-1}A^{-1}\begin{pmatrix} 0 & 0 & A_{13}+A_{55} \\ 0 & 0 & 0 \\ A_{13}+A_{55} & 0 & 0 \end{pmatrix} A^{-1/2}\frac{\partial}{\partial z} \quad (24)$$

$$-\frac{1}{r}\rho^{-1}A^{-1}\begin{pmatrix} A_{11} & 0 & 0 \\ 0 & A_{66} & 0 \\ 0 & 0 & A_{55} \end{pmatrix} A^{-1/2}$$

Due to the diagonal structure of $A^{-1/2}$ the solution to (24) is easily evaluated:

$$P_0 = -\begin{pmatrix} \dfrac{1}{2r} & \dfrac{A_{12}+A_{66}}{A_{11}^{1/2}A_{66}^{1/2}+A_{11}}\dfrac{1}{r}\dfrac{\partial}{\partial\theta} & \dfrac{A_{13}+A_{55}}{A_{11}^{1/2}A_{55}^{1/2}+A_{11}}\dfrac{\partial}{\partial z} \\ \dfrac{A_{12}+A_{66}}{A_{11}^{1/2}A_{66}^{1/2}+A_{66}}\dfrac{1}{r}\dfrac{\partial}{\partial\theta} & \dfrac{1}{2r} & 0 \\ \dfrac{A_{13}+A_{55}}{A_{11}^{1/2}A_{55}^{1/2}+A_{55}}\dfrac{\partial}{\partial z} & 0 & \dfrac{1}{2r} \end{pmatrix}.$$

Therefore, using form (22), we have the following TTBC at $r = const$:

$$\frac{\partial U}{\partial t} + \frac{1}{\sqrt{\rho}}\begin{pmatrix} A_{11}^{1/2}\left(\dfrac{1}{2r}+\dfrac{\partial}{\partial r}\right) & \dfrac{A_{12}+A_{66}}{A_{11}^{1/2}+A_{66}^{1/2}}\dfrac{1}{r}\dfrac{\partial}{\partial\theta} & \dfrac{A_{13}+A_{55}}{A_{11}^{1/2}+A_{55}^{1/2}}\dfrac{\partial}{\partial z} \\ \dfrac{A_{12}+A_{66}}{A_{11}^{1/2}+A_{66}^{1/2}}\dfrac{1}{r}\dfrac{\partial}{\partial\theta} & A_{66}^{1/2}\left(\dfrac{1}{2r}+\dfrac{\partial}{\partial r}\right) & 0 \\ \dfrac{A_{13}+A_{55}}{A_{11}^{1/2}+A_{55}^{1/2}}\dfrac{\partial}{\partial z} & 0 & A_{55}^{1/2}\left(\dfrac{1}{2r}+\dfrac{\partial}{\partial r}\right) \end{pmatrix} U = 0. \quad (25)$$

For VTI case we can take into account the constrain $A_{66} = (A_{11}-A_{12})/2$. Note that other three VTI constrains ($A_{22}=A_{11}$, $A_{23}=A_{13}$, $A_{55}=A_{44}$) do not influence the condition (25).

## 3.2 Three-dimensional Biot's poroelasticity equations in Cartesian coordinates

We consider how to derive TTBCs at a plane boundary normal to one of the coordinate axes $x_i$, $i=1,2,3$. Without viscose terms (i.e. at $\eta=0$) the equations

(8.34) in [2] of waves propagation in anisotropic media in Cartesian coordinates ($i, j = 1, 2, 3$) read

$$\sum_{j=1}^{3} \frac{\partial}{\partial x_j} \left( \sum_{\mu=1}^{3} \sum_{\nu=1}^{3} A_{ij}^{\mu\nu} e_{\mu\nu} + M_{ij} \zeta \right) = \frac{\partial^2}{\partial t^2} \left( \rho u_i + \rho_f w_i \right)$$

$$-\frac{\partial}{\partial x_i} \left( \sum_{\mu=1}^{3} \sum_{\nu=1}^{3} M_{\mu\nu} e_{\mu\nu} + M \zeta \right) = \frac{\partial^2}{\partial t^2} \left( \rho_f u_i + \sum_{j=1}^{3} m_{ij} w_j \right)$$

(26)

Here

$(u_1, u_2, u_3)$ and $\mathbf{w} = (w_1, w_2, w_3)$ are the vectors of solid matrix displacements and relative flow components of the fluid, respectively;

$\rho$ and $\rho_f$ are the mass densities of the porous material and fluid, respectively;

$m_{ij}$ are the so-called efficient densities;

$e_{ij}$ are components of the conventional strain tensor of the solid; $\zeta = -\text{div}\mathbf{w}$ ;

$A_{ij}^{\mu\nu}$, $M_{ij}$, $M$ are the coefficients in the stress-strain matrix relations.

Suppose that all parameters do not depend on time, and the coefficients $A_{ij}^{\mu\nu}$, $M_{ij}$, and $M$ are continuously differentiable. Let $n$ be a one of the coordinates $(x_1, x_2, x_3)$, and $\tau = (\tau_1, \tau_2)$ be the vector of two remained ones. After some obvious algebra, the system (26) can be written in the form

$$-J \frac{\partial^2 U}{\partial t^2} + C^{nn} \frac{\partial^2 U}{\partial n^2} + (C^{n\tau} \nabla_\tau) \frac{\partial U}{\partial n} + C^n \frac{\partial U}{\partial n} + \ldots = 0 \qquad (27)$$

where

$U = (u_1, u_2, u_3, w_1, w_2, w_3)^\mathrm{T}$ is the unknown vector;

$J = \begin{pmatrix} \rho I & \rho_f I \\ \rho_f I & M_f \end{pmatrix}$ is the $6 \times 6$ matrix composed with the help of $3 \times 3$ matrices $M_f = \{m_{ij}\}$ and $I$, the unity matrix;

$C^{nn}$, $C^{n\tau}$, and $C^n$ are the $6 \times 6$ matrices with entries obtained from parameters $A_{ij}^{\mu\nu}$, $M_{ij}$, $M$ and their derivatives by collecting equation coefficients at terms $\dfrac{\partial^2}{\partial n^2}$, $(\cdot \nabla_\tau)\dfrac{\partial}{\partial n}$, and $\dfrac{\partial}{\partial n}$, respectively;

"..." are the remained terms not used while generating TTBCs.

We also suppose that (27) is a hyperbolic system, what is in accordance with the physics of the considered model describing waves' propagation.

Equation (27) differs from (2) by the matrix $J$ at the time derivative. Nevertheless, the theory developed in Section 1 is also applicable to (27). After some algebra we have the following expression for operator $P_1$:

$$P_1 = -(J^{-1} C^{nn})^{-1/2},$$

and equation for operator $P_0$:

$$(J^{-1} C^{nn})^{-1/2} P_0 + P_0 (J^{-1} C^{nn})^{-1/2} =$$
$$= -(C^{nn})^{-1} \left( (C^{n\tau} \nabla_\tau (C^{nn})^{-1/2} J^{1/2}) + (C^{n\tau} (C^{nn})^{-1/2} J^{1/2} \nabla_\tau) + C^n (C^{nn})^{-1/2} J^{1/2} \right).$$

The latter is straightforwardly solved according to formulas (19) – (21). Thus at each boundary point we can evaluate (generally speaking numerically rather than analytically) set of $6 \times 6$ matrices for TTBCs in the form (17) or (22).

## 4 Conclusion

We presented theory of deriving equations for computation of truncated TBCs for second order hyperbolic systems. A way of solving these equations in general case of non-uniform coefficients is developed. The obtained TTBCs are the local (differential) part of the entire TBC operator. Examples of TTBCs for 3D cylindrically anisotropic elasticity and for inviscid poroelasticity (Biot equations) are considered.